\documentstyle[a4,12pt,english]{article} 
\input amssym.def 
\input amssym.tex

\def \sta{{\widetilde{\cal S}} (A)}

\def \A{{\cal A}}
\def \E{{\cal E}}
\def \F{{\cal F}}
\def \H{{\cal H}}
\def \K{{\cal K}}
\def \L{{\cal L}}
\def \M{{\cal M}}
\def \C{{\Bbb C}}
\def \N{{\Bbb N}}
\def \Z{{\Bbb Z}}
\def \R{{\Bbb R}}
\def \S{{\Bbb S}}
\def\sd{\rtimes} 

\def \mapdown#1{\Big\downarrow\rlap{$\vcenter{\hbox
{$\scriptstyle#1$}}$}}

\def \O2{{{\cal O}_2}}

\def \ot{\otimes}

\def \omin{\otimes}
\def \otm{{\mathop{\otimes}\limits^m}_{C(X)}}

\def\cqfd{\vbox{\hrule
height 4pt width 4pt}}
\def\ie{{\it i.e.\/}\ } 
 \def\cf{{\it cf.\/}\ }

\def \rem {\medskip\indent {\mbox{{\it Remark.}}} }
\def \pf {\medskip\indent {\mbox{{\it Proof.}}} }

\def \cst {C$\sp*$}

\title{{\bf A FEW REMARKS ON EXACT C(X)-ALGEBRAS}
} 
\author{Etienne Blanchard}
\date{}

\begin{document}
\maketitle

\renewcommand{\abstractname}{}
\begin{abstract} 

We extend in this paper several results 
of E. Kirchberg, S. Wassermann and the author 
dealing with continuous fields of \cst-algebras to the semi-continuous case. 
We provide a new characterisation of separable lower semi-continuity 
\cst-bundles and we present exactness criteria for $C(X)$-algebras and 
unital bisimplifiable Hopf \cst-algebras. 
\end{abstract}



\bigskip\bigskip\centerline{\large\bf 0.\quad INTRODUCTION}

\bigskip 
\indent\indent 
In a paper on continuous fields of \cst-algebras ([9]), 
E. Kirchberg and S. Wassermann have introduced several characterisations of 
the exactness of a continuous field of \cst-algebras $A$ over a Hausdorff 
compact space $X$ through continuity conditions and fibrewise properties 
of the tensor product with another continuous field. 
Notice that it has been proved later on by E. Kirchberg that all those 
conditions are also equivalent in the separable case to the existence of 
a $C(X)$-linear 
isomorphism between $A$ and a subfield of the trivial continuous field 
$\O2\ot C(X)$, where $\O2$ is the Cuntz \cst-algebra 
([4, Theorem A.1]). 

Our purpose in this article is to analyse these continuity properties 
and to make the difference between those of lower semi-continuity 
and those of upper semi-continuity. 
After a few preliminaries on $C(X)$-algebras, we are led in theorem 2.2 
to characterise, up to a $C(X)$-linear isomorphism, 
every separable lower semi-continuous \cst-bundle 
as a separable $C(X)$-subalgebra of the \cst-algebra $\L (\E )$ of bounded 
$C(X)$-linear 
operators acting on a Hilbert $C(X)$-module $\E$ which admit an adjoint. 
This enables us to give an answer to a question of 
Rieffel on the characterisation of lower semi-continuous 
\cst-bundles in the separable case ([11, page 633]) 
through the generalisation of a previous study of the continuous case 
([3, th\'eor\`eme 3.3]). We then present several criteria of exactness 
for a $C(X)$-algebra which are analogous to those of [9] and we 
construct 
an explicit example of an algebraic tensor product of a continuous field 
of \cst-algebras with a \cst-algebra which cannot be endowed with any 
continuous structure. 
We eventually look at the two notions of exactness for a compact quantum 
group introduced by E. Kirchberg and S. Wassermann, the language of 
multiplicative unitaries allowing us to provide a direct proof of the 
equivalence between these two definitions. 

\medskip 
{\sl I would like to express my gratitude to S Wassermann for fruitful 
discussions held during his different stays at the 
Institut de Math\'ematiques de Luminy. 
I also want to thank S. Baaj for a few comments on a french 
preliminary version.}

\bigskip\bigskip\centerline{\large\bf 1.\quad PRELIMINARIES}

\bigskip
\indent\indent 
Let us recall in this section the basic definitions related to the theory of 
\cst-bundles. 

\medskip\indent {\it Definition} 1.1. ([7]) {\rm 
Given a Hausdorff compact space $X$, a $C(X)$-algebra is a \cst-algebra $A$ 
endowed with a unital morphism from the \cst-algebra $C(X)$ of continuous 
functions on the space $X$ into the centre of the multiplier algebra $\M (A)$ 
of $A$. 

We associate to such an algebra $A$ the unital $C(X)$-algebra 
${\cal A}$ generated by $A$ and $u[C(X)]$ in $\M [A\oplus C(X)]$ 
where $u(g)(a\oplus f)=ga\oplus gf$ for $a\in A$ and $f,g\in C(X)$.} 

\medskip
For $x\in X$, let $C_x(X)$ be the kernel of the evaluation map 
$ev_x :C(X)\rightarrow\C$ at $x$. Denote by $A_x$ the quotient 
of a $C(X)$-algebra $A$ by the closed ideal $C_x(X)A$ and 
by $a_x$ the image of an element $a\in A$ in the fibre $A_x$. 
Then the function 
\begin{center}
$x\mapsto\| a_x\| =\inf\{\| [1-f+f(x)]a\| ,f\in C(X)\}$ 
\end{center}
is upper semi-continuous for all $a\in A$ 
and the $C(X)$-algebra $A$ is said to be a continuous field 
of \cst-algebras over $X$ if the function $x\mapsto\| a_x\|$ is actually 
continuous for every $a\in A$ ([5]). 

\medskip\indent {\it Examples.} 
{\sl 1.} If $A$ is a $C(X)$-algebra and $D$ is a \cst-algebra, 
the spatial tensor product $B=A\ot D$ is naturally endowed with a 
structure of $C(X)$-algebra through the map 
$f\in C(X)\mapsto f\ot 1_{M(D)}\in \M (A\ot D)$. 
In particular, if $A=C(X)$, the tensor product $B$ 
is a trivial continuous field over $X$ with constant fibre $D$. 

\noindent {\sl 2.} If the \cst-algebra $C(X)$ is a von Neumann algebra, 
then any $C(X)$-algebra $A$ is continuous since the lower bound of the 
continuous functions $x\mapsto \| [1-f+f(x)]a\|$, $f\in C(X)$, belongs to 
$C(X)$ for all $a\in A$ in that case. As a consequence, 
if $D$ is a \cst-algebra, the spatial tensor product $A\ot B$ is a continuous 
field over $X$ whose fibre at $x\in X$ is not isomorphic to the tensor product 
$A_x\ot D$ in general. 

\noindent {\sl 3.} Given two $C(X)$-algebras $A_1$ and $A_2$, 
the Hausdorff completion $A_1\otm A_2$ of the algebraic tensor product 
$A_1\odot A_2$ for the semi-norm 
$\|\alpha\|_m=\sup\{\| \sigma_x(\alpha )\| ,\, x\in X\}$, where for every 
$x\in X$, $\sigma_x$ is the map 
$A_1\omin A_2\rightarrow (A_1)_x\omin (A_2)_x$ 
which takes $a\ot b$ to $a_x\ot b_x$   , 
is also a $C(X)$-algebra which satisfies certain minimal properties 
([2, propo\-sition~2.9]). 

\medskip\indent {\it Definition} 1.2. ([3,2]) {\rm
Given a compact Hausdorff space $X$, a $C(X)$-representation of a 
$C(X)$-algebra $D$ in a continuous field of \cst-algebras $A$ over $X$ 
is a $C(X)$-linear morphism $\pi$ from $D$ into 
the multiplier algebra $\M (A)$ of $A$, \ie such that for each $x\in X$, 
the induced representation $\pi_x$ in $\M (A_x)$ factorizes through 
the fibre $D_x$. } 

\medskip 
If the $C(X)$-algebra $D$ admits a $C(X)$-representations $\pi$ in the 
continuous field $A$ over $X$, the function 
\begin{center} 
$x\mapsto \|\pi_x(d)\| =\sup\{\| (\pi (d)a)_x\| , a\in A$ such that 
$\| a\|\leq 1\}$ 
\end{center} 
is lower semi-continuous for all $d\in D$. 

\rem If the induced representation $\pi_x$ of the fibre $D_x$ is faithful 
for each $x\in X$, then the function $x\mapsto\| d_x\| =\|\pi_x(d)\|$ 
is continuous for all $d\in D$ and the $C(X)$-algebra $D$ is therefore 
continuous. 

In particular, a separable $C(X)$-algebra $D$ is continuous if and only if 
there exists a Hilbert $C(X)$-module $\E$ and a $C(X)$-representation $\pi$ 
of $D$ in the continuous field $\K (\E )$ of compact operators acting on $\E$, 
such that for all $x\in X$, the induced representation of the fibre $D_x$ in 
the \cst-algebra $\M (\K (\E )_x)=\L (\E_x)$ (where $\E_x$ is the Hilbert 
space $\E_x=\E\otimes_{ev_x}\C$) is faithful ([3, th\'eor\`eme 3.3]). 

\bigskip\bigskip\centerline{\large\bf 2.\quad 
THE LOWER SEMI-CONTINUITY}

\bigskip 
\indent\indent 
Let us focus on the property of lower semi-continuity 
associated to any $C(X)$-repre\-sen\-tation of a $C(X)$-algebra. 

\medskip
Recall first that a faithful family of 
representations of a \cst-algebra $A$ is a family of stellar 
representations $\{\sigma_\lambda ,\,\lambda\in\Lambda\}$ of $A$ 
such that for any $a\in A$, there exists an index 
$\lambda\in\Lambda$ satisfying $\sigma_\lambda (a)\not =0$. 

\medskip\indent {\it Definition} 2.1. {\rm 
Given a Hausdorff compact space $X$, a lower semi-continuous 
\cst-bundle over $X$ is a pair $(A,\{\sigma_x\} )$ where 
\begin{enumerate}
\item[{\sl a)}] $A$ is a $C(X)$-algebra; 
\item[{\sl b)}] $\{\sigma_x,\, x\in X\}$ is a faithful family of 
representations of the \cst-algebra $A$ such that for every 
$x\in X$, the representation $\sigma_x$ factorizes through the 
fibre $A_x$; 
\item[{\sl c)}] for each $a\in A$, the function $x\mapsto\|\sigma_x(a)\|$ is 
lower semi-continuous. 
\end{enumerate}  }
\setcounter{prop}{1}

\medskip
{\mbox{{\it Remarks.}}} 
\smallskip\noindent {\bf 1.} The faithful $C([0,2])$-representation 
$\sigma$ of the $C([0,2])$-algebra $A=C([0,1])\oplus C([1,2])$ in the Hilbert 
$C([0,2])$-module $\E =C_0([0,2]\backslash\{1\})$ gives us a typical example 
of a separable lower semi-continuous \cst-bundle $(A,\{\sigma_x\} )$ which is 
not continuous 
(we shall come back to this example at the end of the section). 

\smallskip\noindent {\bf 2.} 
If $(A,\{\sigma_x\} )$ is a lower semi-continuous \cst-bundle 
and $\A$ is the associated unital $C(X)$-algebra from definition 1.1, 
the pair $(\A ,\{\widetilde{\sigma}_x\} )$ (where for $x\in X$, 
$\widetilde{\sigma}_x$ is the unital extension of the map $\sigma_x$ 
into the multiplier algebra $\M [\sigma_x(A)\oplus\C ]$) 
is also a unital lower semi-continuous \cst-bundle. 
Indeed, if $\alpha\in\A_+$ is positive, there exist a self-adjoint element 
$b\in A$ and a positive function $f\in C(X)_+$ such that 
$\alpha =b+u[f]$. For $x\in X$, one has the equality 
$Sp(\widetilde{\sigma}_x(\alpha ))\cup\{ 0\}=
[Sp(\sigma_x(b))+f(x)]\cup\{ 0\}$. 
As a consequence, if one writes the decomposition 
$b=b_+-b_-$ where the two positive 
elements $b_+$ and $b_-$ satisfy the relation $b_+b_-=0$, the map 
$x\mapsto\|\widetilde{\sigma}_x(\alpha )\|=\|\sigma_x(b_+)\| +f(x)$ 
is a lower semi-continuous function. 

\bigskip 
We can state the following theorem which generalises a previous 
characterisation of separable continuous fields of \cst-algebras 
(\cf the last remark of the previous section). 

\medskip THEOREM 2.2. {\it
Given a separable lower semi-continuous \cst-bundle 
$(A,\{\sigma_x\} )$ over the Hausdorff compact metrisable space $X$, 
there exist a Hilbert $C(X)$-module $\E$ and a $C(X)$-representation 
$\pi$ of $A$ in the continuous field $\K (\E )$ such that 
for every $x\in X$, one has the isomorphism 
$$\pi_x(A)\simeq\sigma_x(A)\, .$$} 

As in the case of continuous fields, the proof of this theorem 
relies on the following lemma. 

\medskip LEMMA 2.3. {\it 
Let $(A,\{\sigma_x\} )$ be a unital separable lower semi-continuous 
\cst-bundle over the Hausdorff compact space $X$ (whose fibres are all 
assumed to be non zero). Endow the subspace $\sta\subset S(A)$ 
of states $\varphi\in S(A)$ on $A$ whose restriction to the unital 
\cst-subalgebra $C(X)$ is the evaluation map at some point $p(\varphi )\in X$ 
with the restricted weak topology. 

Then the continuous map $p: \sta\rightarrow X$ is open. 
}

\rem The restriction of a state $\varphi\in\sta$ to the ideal 
$\ker (\sigma_x)$, where $x=p(\varphi )\in X$, is always zero. 

\pf Let $\Omega$ be a non empty open set of $\sta$ and consider a point 
$x\in p(\Omega )$. As the \cst-algebra $\sigma_x(A)$ is separable, one can 
find a faithful state $\varphi\in S(\sigma_x(A) )$, positive norm $1$ 
elements $a_2, \cdots a_m$ in $A$ and a constant $\varepsilon\in ]0,1/2[$ 
such that 
\begin{enumerate} 
\item[{\sl a)}] the state $\varphi\circ\sigma_x$ on $A$ belongs to $\Omega$; 
\item[{\sl b)}] the open set 
$V=\{\psi\in\sta ,|(\psi -\varphi\circ\sigma_x )(a_k)|<\varepsilon$ 
for every $2\leq k\leq m\}$ 
is an elementary neighbourhood of the state $\varphi\circ\sigma_x$ in 
$\Omega$. 
\end{enumerate} 
If we set $a_1=1\in A$, then, replacing the indices $2,\ldots ,m$ by a 
permutation if necessary, there exists an integer $n\leq m$ 
such that the family $\{\sigma_x(a_1),\cdots ,\sigma_x(a_n)\}$ is a linearly 
independent family of maximal order. 
Moreover, there exist for every index $j>n$ real numbers $\lambda^j_k$, 
$1\leq k\leq n$, such that the element $b_j=a_j-\sum_{k=1}^{n}\lambda^j_k a_k$ 
satisfies the equality $\sigma_x(b_j)=0$. 
Define the positive sum $a_0=\sum_{j=n+1}^{m}|b_j|\in\ker\sigma_x$ and the 
constants $R=\sup\{ |\lambda^j_k|\}$, 
$\varepsilon_1=\varepsilon /(nR+1)<\varepsilon$. 
Then the open set 
\begin{center} 
$V'=\{\psi\in\sta ,|(\psi -\varphi\circ\sigma_x)(a_k)|<\varepsilon_1$ for all 
$0\leq k\leq n\}$ 
\end{center} 
is an open neighbourhood of $\varphi\circ\sigma_x$ contained in $V$. 
Indeed, if $\psi\in V'$ and $j>n$, then on has 
$|(\psi -\varphi\circ\sigma_x)(a_j)|
\leq\sum_{k=1}^n|\lambda_k^j|\varepsilon_1+\psi (|b_j|)<\varepsilon$. 

\bigskip
Consider for $y\in X$ the linear map $\alpha_y:\C^n\rightarrow \sigma_y(A)$ 
defined by $\lambda =(\lambda_k)\mapsto \sum\lambda_k\sigma_y(a_k)$ 
and let $\S\subset\R^n$ be the the self-adjoint part of the unit sphere 
of $\C^n$ for the norm $\max$. 
By hypothesis, the constant $r=\inf\{\|\alpha_x(\lambda )\| ,\lambda\in\S\}$ 
is strictly positive and the complement $U_1$ of the projection on $X$ 
of the closed set 
$$F_1=\{ (y,\lambda )\in X\times\S , \|\alpha_y(\lambda )\|\leq r/2\}$$ 
is an open neighbourhood of $x$ in $X$ on which the $\sigma_y(a_k)$, 
$1\leq k\leq n$, are still linearly independent. 
As the set 
\begin{center}
$\begin{array}{rcl}
F_2&=&\{ (y,\lambda )\in U_1\times\S,\, \|\alpha_y(\lambda )-n\|\leq n
\,\,{\rm and}\,\, \varphi\circ\alpha_x(\lambda )\leq 0 \} \\
&=&\{ (y,\lambda )\in U_1\times\S,\,\alpha_y(\lambda )\geq 0
\,\,{\rm and}\,\,\varphi\circ\alpha_x(\lambda )\leq 0 \}
\end{array}$ 
\end{center} 
is also closed, we construct in the same way 
an open neighbourhood $U_2\subset U_1$ of $x$ on which the self-adjoint linear 
form $\varphi_y:\alpha_y(\C^n)\rightarrow \C$ defined by the formula 
$\alpha_y(\lambda )\mapsto \varphi\circ\alpha_x(\lambda )$ is a positive form 
of norm $1$ on the operator system $\alpha_y(\C^n)$. 
Let $\varepsilon_2$ be the constant $\varepsilon_2=\varepsilon_1/4$ and 
$U_3$ be the open complement of the projection on $U_2$ of the closed set 
\begin{center}
$F_3=\{ (y,\lambda ,\mu)\in U_2\times\S\times [0, 2n/\varepsilon_2],\, 
\|\alpha_y(\lambda )-\mu\sigma_y(a_0)\|\leq 
(1-\varepsilon_2)\|\alpha_x(\lambda )\|\}$. 
\end{center} 
We are going to construct for each $y\in U_3$ a state 
$\Phi_y\in S(\sigma_y(A))$ such that $\Phi_y\circ\sigma_y\in V'$, 
a result which will end the proof of the lemma. 

\medskip\noindent ---$>$ 
If $\|\sigma_y(a_0)\|\leq\varepsilon_2$, the linear form $\varphi_y$ extends 
by Hahn-Banach to a self-adjoint linear form $\Phi_y$ of norm $\varphi_y(1)=1$ 
on the \cst-algebra $\sigma_y(A)$, whence a state $\Phi_y\circ\sigma_y\in V'$. 

\bigskip\noindent ---$>$ 
If $\|\sigma_y(a_0)\| >\varepsilon_2$, let us first notice that 
$\sigma_y(a_0)\not \in \alpha_y(\R^n)$ : indeed, if there existed a couple 
$(\lambda ,\mu )\in\S\times\R_+^*$ verifying the equality 
$\sigma_y(a_0)=\mu^{-1}\alpha_y(\lambda)$, we would then have the inequality 
$|\mu |>2n/\varepsilon_2$ (since $y\in U_3$), whence the contradiction 
$\|\sigma_y(a_0)\| =\mu^{-1}\|\alpha_y(\lambda)\| 
<\varepsilon_2/2$. 
As a consequence, one can extend the linear form $\varphi_y$ to a self-adjoint 
linear form $\varphi '$ on the operator system $E=\alpha_y(\C^n)+\C a_0$ 
by setting $\varphi_y'(a_0)=0$. 
This unital self-adjoint form satisfies 
the inequality $\|\varphi_y'\|\leq (1-\varepsilon_2)^{-1}$, 
inequality which one only needs 
to check on the self-adjoint part of the unit sphere of $E$.
\begin{itemize}
\item[a)] If $\lambda\in\S$ and $\mu\in [0,1]$, one has the relations 
$\|\alpha_y(\lambda )-\mu\sigma_y(a_0)\|\geq 
(1-\varepsilon_2)\|\alpha_x(\lambda )\|\geq 
(1-\varepsilon_2)|\varphi_y '[\alpha_y(\lambda )-
\mu\sigma_y(a_0)]\|$ because $y\in U_3$; 
\item[b)] if $\lambda\in\S$ and $\nu\in [\varepsilon_2/2n,1]$, then 
$\|\nu\alpha_y(\lambda )-\sigma_y(a_0)\| = 
\nu\|\alpha_y(\lambda )-\nu^{-1}\sigma_y(a_0)\|\geq (1-\varepsilon_2)
|\varphi_y '[\nu\alpha_y(\lambda )-\sigma_y(a_0)] |$ thanks to a); 
\item[c)] if $\lambda\in\S$ and $\nu\in [0,\varepsilon_2/2n[$, the hypothesis 
$\|\nu\alpha_y(\lambda )-\sigma_y(a_0)\| <\|\nu\alpha_y(\lambda )\|$ 
implies the inequality $\|\sigma_y(a_0)\| <2\nu\|\alpha_y(\lambda )\| 
<\varepsilon_2$, which is absurd. 
\end{itemize} 

\medskip 
Choose a self-adjoint extension $\psi_y$ of same norm of the form 
$\varphi_y'$ to the \cst-algebra $\sigma_y(A)$ and write 
the polar decomposition $\psi_y=(\psi_y)_+-(\psi_y)_-$. 
One has by construction the relations 
$(\psi_y)_+(1)+(\psi_y)_-(1)=\|\varphi_y '\|\leq 
(1-\varepsilon_2)^{-1}$ and $(\psi_y)_+(1)-(\psi_y)_-(1)=\psi_y(1)=1$, 
so that $\| (\psi_y)_+-\psi_y\| <\varepsilon_2$ because $(\psi_y)_-(1)\leq 
1/2{\strut {\displaystyle \varepsilon_2}\over
{\strut\displaystyle 1-\varepsilon_2}}\leq\varepsilon_2$ 
(since $\varepsilon_2\leq\varepsilon <1/2$). 
As $\| (\psi_y)_+\|\leq\|\psi_y\|\leq (1-\varepsilon_2)^{-1}<2$, 
the state $\Phi_y =(\| (\psi_y)_+\|)^{-1} (\psi_y)_+$ satisfies 
$\|\Phi_y-\psi_y\|\leq \|\Phi_y-\psi_y\| +\varepsilon_2=
\|\psi_y\| (1-\|\psi_y\|^{-1})+\varepsilon\leq 3\varepsilon_2$ 
and so the state $\Phi_y\circ\sigma_y$ belongs to the neighbourhood $V'$. 
\cqfd
 
\bigskip {\mbox{{\it Proof of 
Theorem 2.2.}}} 
Let us first replace the lower semi-continuous \cst-bundle 
$(A,\{\sigma_x\} )$ by the unital lower semi-continuous \cst-bundle 
$(\A ,\{\widetilde{\sigma}_x\} )$ defined in the remark following 
the definition of lower semi-continuous \cst-bundles, in order to be able 
to apply the previous lemma. 

Given a positive element $a\in \A_+$, a point $x\in X$ and a state 
$\varphi\in S(\sigma_x(\A ))$, corollary 3.7 of [3] enables us to 
construct a continuous section $\Phi$ of the open map 
$p: {\widetilde {\cal S}} (\A )\rightarrow X$ satisfying the equality 
$\Phi_x=\varphi\circ\sigma_x$, whence a $C(X)$-representation $\theta$ of $A$ 
on a Hilbert $C(X)$-module $\F$ such that 
$\varphi\circ\sigma_x(a)\leq\|\theta_x(a)\|\leq\|\sigma_x(a)\|$ 
thanks to proposition 2.13 of [3]. 
Now, an appropriate sum of such $C(X)$-representations allows us 
to construct the desired $C(X)$-representation. 
\cqfd

\bigskip 
One derives from this theorem the following corollary, 
thanks to proposition 4.1 of [2] in particular. 

\medskip COROLLARY 2.4. {\it 
If $(A,\{ \theta_x\} )$ (resp. $(B,\{ \sigma_y\} )$) is a lower semi-continuous 
\cst-bundle over the Hausdorff compact space $X$ (resp. $Y$), 
the family of representations 
$\{ \theta_x\ot \sigma_y\}$ defines a structure of lower semi-continuous 
\cst-bundle with fibres
$\{\theta_x(A)\omin \sigma_y(B)\}$ on the $C(X\times Y)$-algebra $A\omin B$. 

Furthermore, if the representation $\theta_x$ of the fibre $A_x$ is faithful 
for all $x\in X$ and the two topological spaces $X$ and $Y$ coincide, then 
the pair $(A\otm B,\{\theta_x\ot\sigma_x\} )$ is a lower semi-continuous 
\cst-bundle. 
} 

\rem In general, the family $\{\theta_x\ot\sigma_x; x\in X\}$ is not a faithful 
family of representations of the \cst-algebra $A\otm B$ 
(\cf example 3.3.2 of [3]). 

\bigskip\bigskip\centerline{\large\bf 3.\quad 
EXACTNESS CRITERIA FOR A C(X)-ALGEBRA}

\bigskip 
\indent\indent 
In this section, we reformulate the characterisation of exact continuous 
fields obtained by E. Kirchberg and S. Wassermann ([9, theorem 4.6]) 
in the framework of $C(X)$-algebras. 

\medskip PROPOSITION 3.1. {\it 
Given a Hausdorff compact space $X$ and 
a $C(X)$-algebra $A$, the following assertions are equivalent: 
\begin{enumerate}
\item the \cst-algebra $A$ is exact; 
\item for any Hausdorff compact space $Y$ and any $C(Y)$-algebra 
$B$, one has for every couple of points $(x,y)\in X\times Y$ the isomorphism 
$$(A\omin B)_{(x,y)}\simeq A_x\omin B_y\; ;$$ 
\item each fibre $A_x$ is exact and for every \cst-algebra $B$, the 
sequence 
\begin{center}
$0\rightarrow C_x(X)A\omin B\rightarrow A\omin B\rightarrow 
A_x\omin B\rightarrow 0$
\end{center} 
is exact for any point $x\in X$. 
\end{enumerate}

Furthermore, if the space $X$ is metrisable and perfect, 
the previous assertions are equivalent to the following one : 
\begin{enumerate}
\item[2'.] for every $C(X)$-algebra $B$, we have for each point $x\in X$ 
the isomorphism 
\begin{center} 
$(A\omin B)_{(x,x)}\simeq A_x\omin B_x$ 
\end{center} 
\end{enumerate} 
} 

\pf {\em 1}$\Rightarrow${\em 2} 
Consider a $C(Y)$-algebra $B$ and two points $x\in X$, $y\in Y$. 
As the \cst-algebra $A$ is exact, it satisfies property C of Archbold and 
Batty ([12, definition 5.3] or [6]). 
The lines and the columns of the following 
commutative diagram are therefore exact. 
$$\matrix{
&&&0&&&&0&&&&0&&&\cr
&&&\downarrow&&&&\downarrow&&&&\downarrow&&&\cr
0&\!\!\!\rightarrow\!\!\!&\!\! C_x(X)A\!\!\!&\omin&\!\!\!\! C_y(Y)B
\!\!&\!\!\!\rightarrow\!\!\!&\!\! C_x(X)A\!\!\!\!&\omin&\!\! B
\!\!&\!\!\!\rightarrow\!\!\!&C_x(X)A\!\!\!&\omin &\!\!\!\! B_y
&\!\!\!\rightarrow\!\!\!&0\cr
&&&\downarrow&&&&\downarrow&&&&\downarrow&&&\cr
0&\!\!\!\rightarrow\!\!\!&A&\omin&\!\!\!\! C_y(Y)B\!\!
&\!\!\!\rightarrow\!\!\!&A&\omin&\!\! B
\!\!&\!\!\!\rightarrow\!\!\!&A\!\!\!&\omin &\!\!\!\! B_y
&\!\!\!\rightarrow\!\!\!&0\cr
&&&\downarrow&&&&\downarrow&&&&\downarrow&&&\cr
0&\!\!\!\rightarrow\!\!\!&A_x&\omin&\!\!\!\! C_y(Y)B\!\!
&\!\!\!\rightarrow\!\!\!&A_x&\omin&\!\! B\!\!
&\!\!\!\rightarrow\!\!\!
&A_x\!\!\!&\omin &\!\!\!\! B_y&\!\!\!\rightarrow\!\!\!&0\cr
&&&\downarrow&&&&\downarrow&&&&\downarrow&&&\cr
&&&0&&&&0&&&&0&&&
}$$
Consequently, the density of the linear space 
$lin\,\{ C_x(X)\otimes C(Y)+C(X)\otimes C_y(Y)\}$ in the ideal 
$C_{(x,y)}(X\times Y)$ provides us with the exact sequence 
$$0\rightarrow C_{(x,y)}(X\times Y)[A\omin B]\rightarrow A\omin B
\rightarrow A_x\omin B_y\rightarrow 0\; .$$ 

\medskip\noindent 
{\sl 2}$\Rightarrow${\sl 3} Given a $C(Y)$-algebra $B$ and two points 
$x\in X$, $y\in Y$, we have the canonical sequence of epimorphisms 
$(A\omin B)_{(x,y)}\rightarrow (A_x\omin B)_y\rightarrow A_x\omin B_y$, 
whence the exact sequence 
$$ 0\rightarrow A_x\ot C_y(Y)B\rightarrow A_x\ot B\rightarrow 
A_x\ot B_y\rightarrow 0\; .$$ 
In particular, if $Y=\widetilde{\N} =\N\cup\{\infty\}$ is 
the Alexandroff compactification of $\N$ and $B$ is 
the $C(\widetilde{\N} )$-algebra 
$B=\mathop{\prod}\limits_{n=1}^\infty M_n(\C)$, then proposition 4.3 of 
[9] implies the exactness of the fibre $A_x$. 
Now, if the $C(Y)$-algebra $B$ is a trivial continuous field 
$B=C(Y)\otimes D$, 
the sequence of epimorphisms 
$(A\omin B)_{(x,y)}\rightarrow (A\omin D)_x\rightarrow A_x\omin D$ 
provides us with the exact sequence 
$0\rightarrow C_x(X)A\omin D\rightarrow A\omin D\rightarrow 
A_x\omin D\rightarrow 0$. 

\medskip 
Notice that if there exists a sequence of points $x_n\in X\backslash\{ x\}$ 
converging toward the point $x$, then the topological space 
$\{ x_n\}_n\cup\{ x\}\subset X$ is isomorphic to $\widetilde{\N}$, 
whence the implication {\sl 2'}$\Rightarrow${\sl 3}. 

\medskip\noindent 
{\sl 3}$\Rightarrow${\sl 1} Let $B$ be a \cst-algebra, $K$ a two sided closed 
ideal in $B$ and assume that the operator $d\in A\omin B$ belongs to the 
kernel of the quotient map $A\omin B\rightarrow A\omin (B/K)$. 
For all $x\in X$, one then has the relation 
$$ \begin{array}{rl}
d_x\in \ker\{ (A\omin B)_x\rightarrow [A\omin (B/K)]_x\}\!\!\! &=
\ker\{ A_x\omin B\rightarrow A_x\omin (B/K)\} \\ 
&=A_x\omin K\quad {\rm (since}\; A_x\;{\rm is\; exact)}\\ 
&=(A\omin K)_x
\end{array}$$ 
and so $d\in A\omin K$ since the map 
$A\omin B\rightarrow\prod_{x\in X}(A\omin B)_x$ is a monomorphism. 
This means that the sequence 
$0\rightarrow A\omin K\rightarrow A\omin B\rightarrow A\omin (B/K)
\rightarrow 0$ is exact. \cqfd  

\bigskip\noindent 
{\Large {\bf A counter-example}}

\medskip 
Following ideas of [9] and [3], let us construct explicitly 
a separable continuous field of \cst-algebras $A$ on the 
Hausdorff compactification $\widetilde{\N}=\N\cup\{\infty\}$ of $\N$ 
and a \cst-algebra $B$ such that there is no \cst-norm on the algebraic 
tensor product $A\odot B$ which endows this $C(\widetilde{\N})$-module 
with a structure of continuous field. 

\medskip 
Assume that $\Gamma$ is an infinite countable residually finite group 
satisfying property $T$ (for instance $\Gamma =SL_3(\Z)$). 
Consider the countable infinite set of classes of finite dimensional 
irreducible representations $\{\pi_n\}_{n\in\N}$ of $\Gamma$ and 
make the hypothesis that $\pi_0$ is the trivial representation. 
As explained in [9, lemma 4.1], one can then find a strictly 
growing sequence of integers $k_n$, $n\in\N$, with $k_0=0$ , $k_1=1$, 
such that if on sets 
$\sigma_n =\oplus_{i=k_n}^{k_{n+1}-1}\pi_i$, then the limit 
$\lim_{n\rightarrow\infty}\|\sigma_n(a)\|$ exists for all 
$a\in A=(\oplus\pi_n)(C^*(\Gamma ) )$, whence a structure of continuous field 
over $\widetilde{\N}$ for the \cst-algebra $A$ whose fibre at $n\in\N$ is 
the finite dimensional \cst-algebra $A_n=\sigma_n(C^*(\Gamma ) )$. 
Define also the \cst-algebra 
$B=(\oplus_{n\in 2\N}\overline{\sigma}_n)(C^*(\Gamma ) )$, where 
$\overline{\sigma}_n$ denotes the contragredient representation of the 
representation $\sigma_n$ of $\Gamma$. 

If $p\in A$ is the unique projection satisfying the relation $\| p_n\| =1$ 
if and only if $n=0$, the coproduct 
$\delta :C^*(\Gamma )\rightarrow C^*(\Gamma )\omin C^*(\Gamma )$ 
enables us to construct (as in [3, example 3.3.1]) 
the image $q$ of 
the projection $\delta (p)$ in the $C(\widetilde{\N} )$-algebra $A\omin B$. 
For $n$ finite, this projection satisfies the equalities $\| q_n\| =1$ if $n$ 
is even and $\| q_n\| =0$ if $n$ is odd. Consequently, the sequence 
$n\mapsto\| q_n\|$ admits no limits as $n$ goes to $\infty$. 

\bigskip\bigskip\centerline{\large\bf 4.\quad 
UNITAL EXACT HOPF \cst-ALGEBRAS}

\bigskip  
\indent\indent 
We study in this section a quantised presentation of the equivalence 
between the exactness of the reduced \cst-algebra of a discrete group 
$\Gamma$ and the exactness of the group $\Gamma$ obtained by E. Kirchberg and 
S. Wassermann. All the basic definitions and notations of the theory of 
multiplicative unitaries which we shall use may be found in~[1]. 

\bigskip 
Let $(S,\delta )$ be a unital bisimplifiable Hopf \cst-algebra, \ie 
a unital \cst-algebra $S$ endowed with a unital coassociative morphism 
$\delta :S\rightarrow S\omin S$, called coproduct, such that the two 
linear subspaces generated by $\delta (S)(1\otimes S)$ and 
$\delta (S)(S\otimes 1)$ are both dense in the spatial tensor product 
$S\ot S$ (for instance, 
the reduced \cst-algebra $S=C^*_r(\Gamma )$ of a discrete group $\Gamma$). 
Then there exists a Haar measure on $S$ ([14, 13]), 
\ie a state $\varphi$ 
on $S$ satisfying the equalities 
$$\forall a\in S, \quad (\varphi\ot id)\circ\delta (a)=
(id\ot\varphi )\circ\delta (a)=\varphi (a)1\in S.$$ 

Let $(\H_\varphi ,L,e)$ be the G.N.S. construction associated to this state 
$\varphi$ and define the multiplicative unitary 
$\widehat{V}\in L(S)\ot\L (\H_\varphi )\subset 
\L (\H_\varphi\ot\H_\varphi )$, by the formula 
$$\forall a,b\in S,\quad \widehat{V}^*(L(a)e\ot L(b)e)= 
(L\ot L)(\delta (b)(a\ot 1) )(e\ot e).$$ 
The coproduct of the Hopf \cst-algebra $L(S)=\overline{lin}\, 
\{ L(\omega )=(id\ot\omega )(\widehat{V} ),\, \omega\in\L (\H_\varphi )_*\}$ 
satisfies the equality 
$(L\ot L)\delta (a)=\widehat{V}^*(1\ot L(a))\widehat{V}$ for all $a\in S$. 
For $a\in S$, the operator $\lambda (\varphi a)=
(\varphi\ot id)((a\ot 1)\widehat{V} )\in\L (\H_\varphi )$ satisfies the 
relation $1\ot\lambda (\varphi a)=
(\varphi\ot L\ot id)(\delta\ot id)((a\ot 1)\widehat{V} )=
(\varphi\ot L\ot id)(\delta (a)_{12}\widehat{V}_{13})\widehat{V}$ 
since $\varphi$ is a Haar state on $S$. 
As the closed linear span 
$\lambda (\widehat{S})=\{ \lambda (\varphi a), a\in S\}$ defines a non 
degenerate \cst-algebra, the right simplifiability of the Hopf 
\cst-algebra $(S,\delta )$ implies that the linear span 
$lin\, (1\ot\lambda (\widehat{S} ))\widehat{V}^*(L(S)\ot 1)$ is dense in 
$L(S)\ot\lambda (\widehat{S})$ and so the multiplicative unitary 
$\widehat{V}$ is regular ([14]). 
Thus, the \cst-algebra $\lambda (\widehat{S} )$, which is also the 
closed linear span of the elements 
$\lambda (\omega )=(\omega\ot id)(\widehat{V} )$, 
$\omega\in\L (\H_\varphi )_*$, is endowed with a structure of 
bisimplifiable Hopf \cst-algebra for the coproduct 
$(\lambda\ot\lambda )\widehat{\delta}(d)=\widehat{V}(\lambda (d)\ot 1)
\widehat{V}^*$. The couple $(\widehat{S} ,\widehat{\delta} )$ 
will be called the dual Hopf \cst-algebra of the Hopf \cst-algebra 
$(L(S),(L\ot L)\circ\delta )$ ([1]). 

\bigskip
If a \cst-algebra $A$ is endowed with a non-degenerate coaction 
$\delta_A:A\rightarrow\M (A\ot\widehat{S} )$ of the Hopf \cst-algebra 
$\widehat{S}$, a covariant representation of the pair 
$(A,\delta_A)$ is a couple $(\pi ,X)$ where $\pi$ is a representation of $A$ 
on a Hilbert space $\H$ and $X$ is a representation of the Hopf \cst-algebra 
$\widehat{S}$ on the 
same Hilbert space $\H$, \ie a unitary satisfying the relation 
$X_{12}X_{13}\widehat{V}_{23}=\widehat{V}_{23}X_{12}$ in 
$\L(\H\ot\H_\varphi\ot\H_\varphi )$ and such that for all $a\in A$, 
one has the equality $(\pi\ot id)\delta_A(a)=X(\pi(a)\ot 1)X^*$ 
in $\L (\H\ot\widehat{S} )$ (\cf [3, section 5.2]). 
The crossed product $A\sd S$ is then the closed vector space generated 
by the products $\pi (a)(id\ot\omega)(X)\in\L (\H )$, $a\in A$ and 
$\omega\in\L (\H_\varphi )_*$; it is endowed with two canonical 
non-degenerate morphisms of \cst-algebras $\pi :A\rightarrow A\sd S$ 
and $L_X: S_f\rightarrow\M (A\sd S)$, where $S_f$ denotes the full Hopf 
\cst-algebra associated to $S$. 

\medskip 
If one keeps these notations, then one can state the following proposition. 

\medskip PROPOSITION 4.1. {\it  
Given a unital bisimplifiable Hopf \cst-algebra $(S,\delta )$ whose 
dual is the Hopf \cst-algebra $(\widehat{S} ,\widehat{\delta} )$, 
the following assertions are equivalent : 
\begin{enumerate} 
\item the reduced \cst-algebra $L(S)$ is exact; 
\item for every $\widehat{S}$-equivariant sequence 
$0\rightarrow J\rightarrow A\rightarrow A/J\rightarrow 0$, 
the sequence of reduced crossed products 
$$0\rightarrow J\sd_rS\rightarrow A\sd_rS\rightarrow (A/J)\sd_rS
\rightarrow 0$$ 
is also exact. 
\end{enumerate}
} 
 
\pf As the implication $2)\Rightarrow 1)$ is clear (it suffices to look 
at trivial coactions), the main point of the proposition is contained in 
the converse implication. 

\medskip 
Consider a $\widehat{S}$-equivariant sequence 
$0\rightarrow J\rightarrow A\rightarrow A/J\rightarrow 0$ and denote 
by $A\sd_f S$ (resp $J\sd_f S$, $(A/J)\sd_f S$) the full crossed product 
of $A$ (resp. $J$, $A/J$) by $S$. 
Let $\pi$ be the canonical representation of $A$ in 
$A\sd_f S$ and let $X\in\L (A\sd_p S\ot\H_\varphi )$ be the unitary 
defining the representation of $S_f$ in $\M (A\sd_f S)$. 
The quotient $D={(A\sd_f S)/(J\sd_f S)}$ is generated by the products of 
images of elements of $A$ and $S$ in $D$. Furthermore, the image of $J$ in 
$\M (D)$ is zero by construction and so one has an isomorphism 
$D\simeq (A/J)\sd_f S$. 

\medskip 
Notice that the reduced crossed product $A\sd_r S$ is linearly generated 
by the products $\delta_A(a)(id\ot id\ot\omega )(\widehat{V}_{23})\in
\L (A\ot\H_\varphi )$, $a\in A$ and $\omega\in\L (\H_\varphi )_*$ 
([1]) and that the unitary $\widehat{V}$ belongs to the \cst-algebra 
$S\ot\L (\H_\varphi )$. As a consequence, one defines a monomorphism 
$\Phi_A :A\sd_r S\rightarrow A\sd_f S\ot S$ through the formula 
$$\matrix{
\Phi_A\left[ \delta_A(a)(id\ot id\ot\omega )(\widehat{V}_{23})\right] 
&=X^*(\pi\ot id)\left[ \delta_A(a)(id\ot id\ot\omega )(\widehat{V}_{23})
\right] X\hfill \cr
&=(\pi (a)\ot 1)(id\ot id\ot\omega )(X_{13}\widehat{V}_{23})
\in A\sd_pS\ot S\, .\hfill}$$ 
One constructs in the same way monomorphisms 
$\Phi_J :J\sd_r S\rightarrow J\sd_f S\ot S$ and 
$\Phi_{A/J} :(A/J)\sd_r S\rightarrow (A/J)\sd_f S\ot S$. 

As the \cst-algebra $L(S)$ is by assumption exact, one has the following 
commutative diagram whose lower line is exact. 
$$\matrix{
&0&&0&&0&\cr
&\downarrow &&\downarrow &&\downarrow &\cr 
0\; \rightarrow &J\sd_r S&\rightarrow &A\sd_r S&\rightarrow &(A/J)\sd_r S& 
\rightarrow\; 0\cr
&\mapdown{\Phi_J} &&\mapdown{\Phi_A} &&\mapdown{\Phi_{A/J}} &\cr
0\; \rightarrow &J\sd_f S\ot L(S)&\rightarrow &A\sd_f S\ot L(S)&\rightarrow 
&(A/J)\sd_f S\ot L(S)& \rightarrow\; 0 
}$$

Assume that $T\in A\sd_rS$ belongs to the kernel 
$\ker (A\sd_rS\rightarrow (A/J)\sd_r S)$ and let $(u_\lambda )_\lambda$ be 
an approximate identity of the ideal $J$ in $A$. 
Then the net $\|\Phi_A (T)[1-(\pi (u_\lambda )\ot 1)]\| =
\|\Phi_A (T[1-\delta_A (u_\lambda )])\|$ converges to zero and so 
$T\in J\sd_rS$, whence the expected exactness of the upper line. 
\cqfd 

\medskip 
As noticed in the latest remark of [9], one deduces from this 
proposition the following corollary which generalises corollary 5.10 
of [3], where the amenable case was considered. 

\medskip CORROLARY 4.2. {\it  
If $S$ is a unital bisimplifiable Hopf \cst-algebra whose reduced 
\cst-algebra $L(S)$ is exact and $A$ is a $C(X)$-algebra endowed with 
a non-degenerate $C(X)$-linear coaction 
$\delta_A:A\rightarrow \M (A\otimes\widehat{S})$ 
of the (discrete) quantum dual $\widehat{S}$ of $S$, 
then the fibre at $x\in X$ of the reduced crossed product 
$A\sd_rS$ is $(A\sd_rS)_x=A_x\sd_rS$. 

In particular, if the $C(X)$-algebra $A$ is continuous, 
the $C(X)$-algebra $A\sd_rS$ is a continuous field of \cst-algebras with 
fibres $A_x\sd_rS$. 
}

\rem It would be interesting to know whether proposition 4.1 
generalises to the framework of Hopf \cst-algebras $S_V$ associated 
to a regular multiplicative unitary $V\in\M (\widehat{S}_V\ot S_V)$, \ie 
given a $\widehat{S}_V$-equivariant sequence 
$0\rightarrow J\rightarrow A\rightarrow A/J\rightarrow 0$, 
does the exactness of the \cst-algebra $S_V$ imply the exactness of the 
sequence of reduced crossed products $0\rightarrow J\sd_rS_V\rightarrow 
A\sd_rS_V\rightarrow (A/J)\sd_rS_V\rightarrow 0$ ?

\renewcommand{\refname}{

\bigskip
\small {\it Institut de Math\'ematiques de Luminy\\
CNRS--Luminy case 907\\ 
F--13288 Marseille CEDEX 9\\ 
\\E-mail: E.Blanchard@iml.univ-mrs.fr

}
\end{document}